\documentstyle{jsc}

\newcommand{\NN}{{\bf N}}
\newcommand{\T}{{\bf T}}
\newcommand{\Z}{{\bf Z}}
\newcommand{\gr}{Gr\"obner }
\newcommand{\sa}{subalgebra }

\newcommand{\k}{k[X]}
\newcommand{\r}{R[X]}
\newcommand{\ry}{R[Y]}
\newcommand{\<}{\langle}
\def\>{\rangle}
\newcommand{\lra}{\longrightarrow}

\newcommand{\Xb}{X^{\vec{\beta}}}
\newcommand{\redF}{\stackrel{F}{\longrightarrow}}
\newcommand{\redH}{\stackrel{H}{\longrightarrow}}
\newcommand{\sig}{\stackrel{G}{\lra}_{\mbox{\scriptsize{\rm si}}}}
\newcommand{\va}{{\vec{\mbox{\scriptsize{\rm a}}}}}
\newcommand{\Va}{{\vec{\mbox{\rm a}}}}
\newcommand{\Vh}{{\vec{\mbox{\rm h}}}}
\newcommand{\Vs}{{\vec{\mbox{\rm s}}}}
\newcommand{\VQ}{{\vec{Q}}}

\newcommand{\VP}{{\vec{P}}}
\newcommand{\vpi}{{\vec{\pi}}}
\newcommand{\val}{{\vec{\alpha}}}
\newcommand{\veta}{{\vec{\eta}}}

\newcommand{\e}{{\vec{\mbox{\scriptsize{\rm e}}}}}
\newcommand{\Ve}{{\vec{\mbox{\rm e}}}}
\newcommand{\Vapr}{{\vec{\mbox{\rm a}}}\hspace{1pt}'}
\newcommand{\VQpr}{{\vec{\mbox{\rm Q}}}\hspace{1pt}'}
\newcommand{\Vhpr}{{\vec{\mbox{\rm h}}}\hspace{1pt}'}

\newcommand{\het}{\mbox{\rm ht}}
\newcommand{\lt}{\mbox{\rm lt}}
\newcommand{\vlt}{{\vec{\mbox{\rm lt}}}}
\newcommand{\lp}{\mbox{\rm lp}}

\newcommand{\sLp}{\mbox{\scriptsize{\rm Lp}}}

\newcommand{\lc}{\mbox{\rm lc}}
\newcommand{\Lt}{\mbox{\rm Lt}}
\newcommand{\Lp}{\mbox{\rm Lp}}
\newcommand{\Lc}{\mbox{\rm Lc}}

\newcommand{\syz}{\mbox{\rm Syz}}
\newcommand{\lts}{\mbox{\rm LtSyz}_A}

\newtheorem{th}{Theorem}[section]
\newtheorem{deftn}[th]{Definition}
\newenvironment{df}{\begin{deftn} \rm}{\end{deftn}}
\newtheorem{exx}[th]{Example}
\newenvironment{ex}{\begin{exx} \rm}{$\triangle$ \end{exx}}

\newtheorem{cor}[th]{Corollary}
\newtheorem{prop}[th]{Proposition}
\newenvironment{pf}{{\bf Proof.} }{$\Box$ \vspace{.2in}}

\newtheorem{algo}{Algorithm}
\newenvironment{mycaption}[1]{ \begin{algo}\begin{center}:  
\rm #1 \end{center} }{\end{algo}}

\title [Analogs of \gr Bases in Polynomial Rings over a Ring]
{Analogs of \gr Bases in Polynomial Rings over a Ring}
\author[J. L. Miller]{J. Lyn Miller\\Western Kentucky University}

\begin{document}

\maketitle

%
%
%
%
%
%
\begin{abstract}
In this paper we will define analogs of \gr bases for
$R$-subalgebras and their ideals in a polynomial ring $R[x_1,\ldots,x_n]$
where $R$ is a noetherian 
integral domain with multiplicative identity and in which
we can determine ideal membership and compute syzygies.
The main goal is to present and verify algorithms for constructing
these \gr basis counterparts.
As an application, we will produce a method for
computing generators for the first syzygy module of a
subset of an $R$-\sa of $R[x_1,\ldots,x_n]$ where each coordinate
of each syzygy must be an element of the subalgebra.
\end{abstract}

%
%
%
%
\section{Introduction}
The concept of \gr bases for ideals of a polynomial ring over a field $k$
can be adapted in a natural way for $k$-subalgebras of such a polynomial
ring.  Robbiano and Sweedler (refer to \cite{RS}; see also \cite{KM}) defined a
SAGBI basis\footnote{The name SAGBI is
an acronym standing for Subalgebra Analog to \gr Bases for Ideals.}
for a $k$-subalgebra $A$ of $k[x_1,\ldots,x_n]$ to be a subset $F\subseteq A$
whose leading power products generate the multiplicative monoid of
leading power products of $A$.  
The properties and applications of SAGBI bases strongly imitate
many of the standard \gr basis results when a suitable accompanying 
reduction algorithm is defined.  Sweedler (see \cite{IVR}) went on
to extend the theory of \gr bases in a way that can be used to define
them for ideals of $k$-subalgebras of $k[x_1,\ldots,x_n]$; this was
briefly presented more explicitly by Ollivier (see \cite{Oll}).
Based on their work, we define a SAGBI-\gr basis for an ideal $I$
of a $k$-\sa $A\subseteq k[x_1,\ldots,x_n]$ to be a subset $G\subseteq I$
whose leading power products generate the monoid-ideal consisting
of the leading power products of $I$ in the monoid of those of $A$.
Basic properties and applications of SAGBI-\gr bases are again
straight-forward adaptations of the usual \gr basis theory.
(See also \cite{Mil}.)

Our aim in this paper is to extend the theory of SAGBI and SAGBI-\gr
bases to the context of a polynomial ring over a 
noetherian integral domain $R$ in which we can determine ideal
membership and compute syzygies.
As we know from the study of this same extension process for \gr bases,
the leading coefficients of the polynomials now play a large role.
The definitions, results, and especially techniques in this new setting
are no longer such carbon copies of those for \gr bases,
although we always attempt to parallel them as much as possible.
In particular, the definition of a SAGBI basis in $R[x_1,\ldots,x_n]$
must now allow for addition of leading terms, not just multiplication.
Therefore, the monoid of leading power products used for SAGBI bases
in $k[x_1,\ldots,x_n]$ must be exchanged for a much larger structure,
namely, the $R$-\sa that it generates in $R[x_1,\ldots,x_n]$.
Likewise, for SAGBI-\gr bases in $R[x_1,\ldots,x_n]$, the monoid-ideal
in the definition over $k[x_1,\ldots,x_n]$ must be enlarged to an
ideal of the new $R$-subalgebra just mentioned.

The main goals of this paper are to present and verify algorithms
for constructing SAGBI and SAGBI-\gr bases in $R[x_1,\ldots,x_n]$,
as well as outlining some of their basic properties.
As an application, we will also present a method for
computing generators for the first syzygy module of a
subset of an $R$-\sa of $R[x_1,\ldots,x_n]$ where each coordinate
of each syzygy must be an element of the subalgebra.

\section{Notation}
Our context is the polynomial ring $R[x_1,\ldots,x_n]$ in $n$ variables,
where $R$ is a noetherian integral domain in which we can determine ideal
membership and compute syzygies.  (When the coefficient ring is a field,
we use the symbol $k$ instead of $R$.)
We abbreviate this polynomial ring
as $\r$.
The notation $R[S]$ stands for the $R$-subalgebra generated by
the subset $S\subseteq \r$.
Throughout this paper, $A$ is an $R$-\sa of $\r$.

The symbol $\NN$ represents the non-negative integers, and $\T_X$
represents the set of all power products $\prod_{i=1}^{n} x_{i}^{\beta_i}$
with $\beta_i\in \NN$ of the variables $x_1,\ldots, x_n$.
We will often
abbreviate such a power product as $X^{\vec{\beta}}$
where $\vec{\beta}$ is the exponent vector $(\beta_1,\ldots,\beta_n)$.
More generally, we have
\begin{df}\label{1}
Let $S\subseteq\r$.
  An $S$-power product is a (finite) product
of the form $s_1^{e_1} \cdots s_m^{e_m}$ where $s_i \in S$
and $e_i \in \NN$ for $1\leq i \leq m$.  We usually 
write this simply as
$S^\e$, where $\Ve$ represents that
vector in $\oplus_S \NN$ whose coordinates are all 0 except for
$e_1, \ldots , e_m$ in the positions corresponding to
$s_1, \ldots , s_m$.
\end{df}

\begin{df}
Given a term order on $\r$,
$p\in\r$, and $S\subseteq\r$, we define
\begin{eqnarray*}
\lp(p) &=& \mbox{the leading $X$-power product of $p$}\\
\lc(p) &=& \mbox{the leading coefficient of $p$}\\
\lt(p) &=& \lc(p)\lp(p) = \mbox{ the leading term of $p$}\\
\Lp S&=&\{\lp(s):s\in S\}
\end{eqnarray*}
while $\Lc S$ and $\Lt S$ are similarly defined.
We also establish the convention that $\lp(0)$ is 
undefined while $\lc(0)$ and $\lt(0)$ are 0.
\end{df}

We borrow the following terminology from \cite{RS}.
\begin{df}\label{315}
Let $S\subseteq\r$.  
Given an expression $\sum_{i=1}^N r_i s_i$ with $r_i\in R$
and $s_i\in S$, we define its height, written 
$\het(\sum_{i=1}^N r_i s_i)$,
to
be $\max_i \lp(s_i)$.
Moreover,
we say that $s_{i_0}$ contributes to the height of the expression
 if $\lp(s_{i_0}) = \max_i \lp(s_i)$.
\end{df}
We emphasize that the height is defined only for specific representations
of an element of $\r$, not for that element itself.
Finally, 
we establish the following notation:
\begin{df}
For an $R$-\sa $A\subseteq\r$ and a subset $S\subseteq A$,
\begin{enumerate}
\item $\<S\>_A$ represents the ideal
of $A$ generated by $S$, omitting the subscript when $A$
is obvious.  
\item $\syz_A(S)=\{\Va=(a_s)_{s\in S}\in \oplus_S A:
\sum_{s\in S}a_s s=0\}$, the $A$-module of syzygies of $S$ whose coordinates
all belong to $A$.  We call an element of $\syz_A(S)$ an $A$-syzygy of $S$.
\item If $A$ is a graded algebra,
and $\deg(a)$ represents the degree of $a\in A$,
then $\syz^*_A(S)=
\{\Va\in \syz_A(S): \deg(a_s s) \mbox{ is the same }\forall a_s s \neq 0\}$.
This common value of $\deg(a_s s)$ is called the degree of the syzygy, and
we write it as $\deg(\Va)$.  The elements of $\syz_A^*(S)$ are called
homogeneous $A$-syzygies of $S$.
\end{enumerate}
The subscripts in $\syz_A(S)$ and $\syz_A^*(S)$ 
will be omitted when $A$ is obvious.
\end{df}

\section{SAGBI Bases in $\r$}
Our first goal is to define a SAGBI
basis and present an algorithm
for its construction.
\begin{df}\label{302}
Let $A$ be an $R$-\sa of $\r$.  
We say that $F\subseteq A$ is a SAGBI basis for $A$ if $\Lt F$ generates
the $R$-\sa $R[\Lt A]$, i.e. $R[\Lt A]= R[\Lt F]$.
\end{df}

We consider an operation which parallels the reduction algorithm used
in \gr basis theory.
\begin{df}\label{307}
Let $g \in \r$, and let $F \subseteq \r$.  We will say that
$g$ s-reduces to $h$ via $F$ in one step,
written $g \redF h$, if there exist a non-zero term $c\Xb$ of $g$ and
$F$-power products $F^{\e_1}, \ldots , F^{\e_N}$
such that
\begin{enumerate}
\item
$\lp(F^{\e_i}) = \Xb$ for $1 \leq i \leq N$.
\item
$c = \sum_{i=1}^N r_i \lc(F^{\e_i})$ where $r_i \in R$ for $1 \leq i \leq N$.
\item
$h=g- \sum_{i=1}^N r_i F^{\e_i}$.
\end{enumerate}
We also write $g \redF h$ if there is a finite chain of 1-step
s-reductions leading from $g$ to $h$;
we say that $g$ s-reduces to $h$ via $F$ in this case.
If $h$ cannot be further s-reduced via $F$, then we call it
a final s-reductum of $g$.
\end{df}
It is obvious that if
$g \redF h$, then $g-h \in R[F]$.
Well-ordering of $\T_X$ implies that any chain of 1-step s-reductions
must terminate.

To s-reduce $g \in \r $ via a finite set $F$ requires us to do two
things at each step.  After we have chosen the
term $c\Xb$ of $g$ that we wish to eliminate, we must be able
to tell
\begin{enumerate}
\item
whether $\Xb$ lies in the multiplicative monoid generated by $\Lp F$, and
\item
whether $c$ belongs to the ideal of $R$ generated
by $\{\lc(F^\e):\lp(F^\e)=\Xb\}$.
\end{enumerate}
To address the first issue, we need to
solve the inhomogeneous linear diophantine system
arising from the exponents of the variables in
$\Xb = \lp(F^{\vec{\epsilon}})$
for $\vec{\epsilon} \in \oplus_F\NN$.\footnote
{Refer to \cite{Dach} for a subroutine that can determine such solutions.}
To address the second point simply requires
that we determine ideal membership in $R$, which we have assumed is possible.

By a standard proof, we can also show
\begin{prop}\label{309}
The following are equivalent for $F\subseteq A$:
\begin{enumerate}
\item $F$ is a SAGBI basis for $A$
\item For every $a \neq 0\in A$, the final s-reductum of $a$ via $F$ is always 0.
\item Every $a\in A$ has a SAGBI representation with respect to $F$,
that is, a representation $$a=\sum_{i=1}^N r_iF^{\e_i},{~~}r_i\in R$$
such that $\max_i\lp(F^{\e_i})=\lp(a)$.
\end{enumerate}
\end{prop}
\begin{cor}\label{311}
A SAGBI basis for $A$ generates $A$ as an $R$-\sa.
\end{cor}
\begin{cor}
Suppose $F$ is a SAGBI basis for $A$.  An element $p\in \r$
belongs to $A$ $\iff$ $p\redF 0$.
\end{cor}

Now we write $A=R[F]$, where $F=\{f_1,f_2,\ldots\}$ is not necessarily
finite.
To design an algorithm for constructing a SAGBI
basis for $A$, we intend
to determine a collection of polynomials
related to $F$ such that if each of these
polynomials s-reduces to 0 via $F$,
then $F$ is a SAGBI basis.  These polynomials
mimic the S-polynomials of ordinary \gr basis theory,
and this desired property
will be
the basis of our construction algorithm.

Represent $A=R[F]$ as
the homomorphic image of a polynomial ring $R[Y]$ (where the cardinality
of $Y=\{y_1,y_2,\ldots\}$ is the same as that of $F$)
via the usual evaluation homomorphism sending each $y_i\mapsto f_i$.
We will now equip $R[Y]$ with a 
graded $R$-module structure (which may not be 
based on any term order in $R[Y]$).
Given
$P(Y) \in \ry$, we define
$$\deg P(Y) = \max \{ \lp(F^\va):Y^\va\mbox{ occurs in }P(Y) \}.\footnote{
It is not necessarily true that $\deg P(Y) = \lp(P(F))$.
For example, if $F=\{f_1,f_2\}=\{x^2,x^2+1\}\subseteq R[x]$,
then $\deg(y_2-y_1)=x^2$, whereas $\lp(f_2-f_1) = 1$.}$$
It is easy to check that this degree
map from $\ry \rightarrow \T_X$ truly does give a
grading on $\ry$.
Notice that the homogeneous elements with respect to
this presumed grading will be those polynomials $P(Y)$
whose terms give rise to $F$-monomials all having the same
leading $X$-power product.  

Now define an evaluation map
$\pi:\ry \rightarrow R[\Lt F]$ via $y_i\mapsto \lt(f_i)$.
The ideal
$I(\Lt F)=\{P(Y):\pi(P(Y))=P(\Lt F)=0\}$
is homogeneous with respect to the 
$\T_X$-grading on $\ry$.
Its homogeneous generators take the
place in our current theory of the usual S-polynomials.
Recall that such generators may be computed using
the familiar tag variable technique of ordinary \gr basis
theory.  (Refer to \cite{AL} et al.)

We are now in a position to state and prove the main
result of this section.

\begin{th}\label{314}
Let $F \subseteq \r$ have distinct elements, and let
$\{ P_j(Y):j \in J \}$ be a set of
$\T_X$-homogeneous generators for $I(\Lt F) \subseteq \ry$.
$F$ is a SAGBI basis for $R[F]$ $\iff$
for each $j \in J$, $P_j(F) \redF 0$.
\end{th}
\begin{pf}
$\Longrightarrow$: This direction is a trivial corollary of Proposition
\ref{309}.

$\Longleftarrow$:
Let $h \in R[F]$.  We will show that $\lt(h) = \sum_i r_i \lt(F^{\e_i})
\in R[\Lt F]$, which will fulfill Definition \ref{302}.

Write $h = \sum_{i=1}^m c_i F^{\e_i}$; 
furthermore, we will assume that
this representation has the smallest possible
height $t_0 = \max_i \lp(F^{\e_i})$
of all such representations.
We know that $\lp(h) \leq t_0$.
Suppose that $\lp(h) < t_0$;
without loss of generality, let the first $N$
summands be the ones for which $\lp(F^{\e_i}) = t_0$.
Then cancellation of their leading $X$-power products
must occur; i.e. $\sum_{i=1}^N c_i \lt(F^{\e_i}) =0$.
Hence, we obtain an element $P(Y)=\sum_{i=1}^N c_i Y^{\e_i} \in I(\Lt F)$.
We can then write 
\begin{equation}\label{319}
\sum_{i=1}^N c_i Y^{\e_i} = P(Y) = \sum_{j=1}^M g_j(Y) P_j(Y)
\end{equation}
where the elements $P_j(Y)$ are the stated generators of $I(\Lt F)$
and the polynomials $g_j(Y)\in R[Y]$.  Moreover, 
we
may assume that each $g_j(Y)$ is $\T_X$-homogeneous 
(since $P(Y)$ and every $P_j(Y)$ are) and also that
$\deg [g_j(Y)P_j(Y)]= \deg P(Y) = t_0$ for $1\leq j \leq M$.

We have assumed that each $P_j(F) \redF 0$; therefore, we
have SAGBI representations 
$P_j(F) = \sum_{k=1}^{n_{kj}} c_{kj} F^{\e_{kj}}$.
By definition, these sums must
have heights $\max_k \lp(F^{\e_{kj}})
= \lp(P_j(F)) < \deg P_j(Y)$ for each $j$, where
the last inequality holds because $P_j(Y)\in 
I(\Lt F)$, so that the highest $X$-terms in $P_j(F)$ cancel.
Then for each $j$, $1\leq j \leq M$,
\begin{equation}\label{318}
g_j(F)P_j(F) = \sum_{k=1}^{n_j} c_{kj} g_j(F) F^{\e_{kj}} 
\end{equation}
Define $t_j$ to be the height of the
right-hand sum in Equation (\ref{318}), and observe that
$$t_j \leq \deg g_j(Y) \cdot \max_k \lp(F^{\e_{kj}})
 <  \deg g_j(Y) \cdot \deg P_j(Y)
 = t_0.$$
Note that
it is impossible
for $F^{\e_1}$ to occur in the right-hand sum
of Equation (\ref{318})
since $t_0 = \lp(F^{\e_1})$.

Returning to our representation 
in Equation (\ref{319}),
we define $d_j$ to be the coefficient of $Y^{\e_1}$
in $g_j(Y)P_j(Y)$ and assume that $d_j \neq 0$
for $1 \leq j \leq M_1$, $d_j =0 $ for $j>M_1$.
Furthermore, let us define $U_j(Y) = g_j(Y)P_j(Y) -d_j Y^{\e_1}$;
we solve this equation
for $d_jY^{\e_1}$, apply the evaluation map $R[Y]\longleftarrow R[F]$,
and substitute using Equation (\ref{319}) to obtain
$$d_j F^{\e_1} = - U_j(F) + \sum_{k=1}^{n_j} c_{kj} g_j(F) F^{\e_{kj}},
{~~~~~}1\leq j \leq M.$$
Observe that $F^{\e_1}$ may not occur on the right-hand
side of the equation: it did not appear on the right-hand-side
of Equation (\ref{319}), and $U_j(Y)$ contains no term involving
$Y^{\e_1}$, whence $U_j(F)$ contains no term involving $F^{\e_1}$
(This last statement requires our assumption that the members of $F$
are distinct.)

Our definition of $d_j$ implies that $c_1 = \sum_{j=1}^M d_j$.
Therefore,
$$c_1 F^{\e_1} = \sum_{j=1}^M d_j F^{\e_1} = \sum_{j=1}^M 
\left[ -U_j(F)  + \sum_{k=1}^{n_j} c_{kj} g_j(F) F^{\e_{kj}} \right].$$
We can now replace $c_1 F^{\e_1}$ by this sum in
the original expression for our polynomial $h$ to get
\begin{eqnarray*}
h &=& \sum_{j=1}^M
\left[ -U_j(F)  + \sum_{k=1}^{n_j} c_{kj} g_j(F) F^{\e_{kj}} \right]
+ \sum_{i=2}^m c_i F^{\e_i}\\
 &=& \sum_{j=1}^M [ -g_j(F)P_j(F) + d_j F^{\e_1} ]
+\sum_{j=1}^M \left[ \sum_{k=1}^{n_j} c_{kj} g_j(F) F^{\e_{kj}} \right]
+ \sum_{i=2}^m c_i F^{\e_i}\\ 
 &=& \sum_{i=1}^N [ -c_i F^{\e_i} ]  + c_1 F^{\e_1}
+\sum_{j=1}^M \left[ \sum_{k=1}^{n_j} c_{kj} g_j(F) F^{\e_{kj}} \right]
+ \sum_{i=2}^m c_i F^{\e_i}\\ 
 &=& \sum_{i=N+1}^m c_i F^{\e_i} 
+\sum_{j=1}^M \left[ \sum_{k=1}^{n_j} c_{kj} g_j(F) F^{\e_{kj}} \right].
\end{eqnarray*}
If we examine this final expression closely, we see that its height
is strictly less than that of our initial representation
for $h$, for
\begin{enumerate}
\item
The height of $\sum_{i=N+1}^m c_i F^{\e_i}$
is strictly less than the old maximum, $t_0$, by
choice of $N$.
\item
We have already seen that the height of
the second sum, which is the maximum of the $t_j$ 
we worked with above, is strictly less than $t_0$.
\end{enumerate}
But this contradicts our initial assumption that we had
chosen a representation for $h$ that had the smallest
possible height.
Thus, $F$
is a SAGBI basis for $R[F]$.
\end{pf}

We may now present an algorithm for
computing SAGBI bases.  See Algorithm \ref{316}.

\begin{figure}[hbt]
\begin{center}
\frame{
\begin{picture}(375,245)
\normalsize{
\put(30,220){{\bf INPUT:} $F$}
\put(30,195){{\bf OUTPUT:} A SAGBI basis for $R[F]$}
\put(30,170){{\bf INITIALIZATION:} $H:=F,{~}oldH:=\emptyset$}
\put(30,145) {{\bf WHILE} $H \neq oldH$ {\bf DO}}
\put(50,120) {$Y := \{y_i:h_i\in H\}$, a set of variables}
\put(50,95) {Choose a $\T_X$-homog. generating set $\cal P$ for $I(\Lt H)$ in $\ry$.}
\put(50,70) {$red{\cal P} := \{\mbox{final s-reducta via $H$ of }P(H) :
              P(Y)\in {\cal P} \}- \{0\}$}
\put(50,45) {$oldH := H$}
\put(50,20) {$H := H \cup red{\cal P}$}
}
\end{picture}
}
\begin{mycaption}{SAGBI Basis Construction Algorithm}\label{316}
\end{mycaption}
\end{center}
\end{figure}

Theoretically, Algorithm \ref{316} can be used with an infinite
input set $F$ because all our results so far have been
carefully designed not to require any finiteness conditions.
Thus, if we assume that we can find generators for
$I(\Lt F)$ when $F$ is infinite (which may be quite a
stretch of imagination!), we shall see that it
makes sense to apply the algorithm to any size input set.
To this end,
we validate that the algorithm
produces a SAGBI basis, although it
{\bf need not terminate}, even with finite input.
(See \cite{RS} for a discussion of infinite
SAGBI bases in $\k$.)

\begin{prop}\label{333}
Let $H_\infty = \cup H$ over all passes of the WHILE loop.
Then $H_\infty$ is a SAGBI basis for $R[F]$.
Moreover, if $F$ is finite and $R[F]$ has a finite SAGBI basis,
then Algorithm \ref{316} terminates and produces a finite SAGBI
basis for $R[F]$.
\end{prop}
\begin{pf}
Set ${\cal P}_\infty = \cup {\cal P}$ over all passes of the loop,
and let $Y_\infty$ be a set of variables $y_i$, one for each element
$h_i\in H_\infty$.
We will show that ${\cal P}_\infty$ is a set of $\T_X$-homogeneous
generators for $I(\Lt H_\infty)\subseteq R[Y_\infty]$,
and then that
each element of ${\cal P}_\infty$ s-reduces to 0 via $H_\infty$.

${\cal P}_\infty$ is $\T_X$-homogeneous
since each $\cal P$ of each loop is.
Now choose $P(Y_\infty) \in I(\Lt H_\infty)$.  Since
only finitely many $y_i$
can occur in $P(Y_\infty)$, only finitely many
$h_i \in H_\infty$ occur in $P(H_\infty)$.
The sets $H$ are nested, so
these particular $h_i$'s must all belong to the set
$H=H_{N_0}$ produced by the end of some finite number $N_0$ of
loops.  Let ${\cal P}_{N_0}$ denote the generating set for $I(\Lt H_{N_0})$,
and let $Y_{N_0}\subseteq Y_\infty$ be the subset of variables
corresponding to $H_{N_0}$.
Then $P(Y_\infty) \in I(\Lt H_{N_0})=
\<{\cal P}_{N_0}\>_{R[Y_{N_0}]}\subseteq 
\<{\cal P}_\infty\>_{R[Y_\infty]}$.
Hence, $I(\Lt  H_\infty) \subseteq \<{\cal P}_\infty\>$.
Conversely, each element $P(Y_\infty)$ of ${\cal P}_\infty$ belongs to the 
set ${\cal P}_{N_0}$
of some pass of the WHILE loop; whence, 
$P(Y_\infty) \in I(\Lt H_{N_0})\subseteq I(\Lt H_\infty)$.
Thus, ${\cal P}_\infty \subseteq I(\Lt H_\infty)$, and
$\< {\cal P}_\infty \> = I(\Lt H_\infty)$.

We have just pointed out that if $P(Y_\infty) \in {\cal P}_\infty$,
then we may assume that $P(Y_\infty) 
\in I(\Lt H_{N_0})$ for some pass, in this case the $N_0$-th, of the loop.
Clearly, either $P(Y_\infty) \stackrel{H_{N_0}}{\lra} 0$
or $P(Y_\infty) \stackrel{H_{N_0+1}}{\lra} 0$.
In either case, we have $P(Y_\infty) \stackrel{H_\infty}{\lra} 0$.
Thus, by
Theorem \ref{314}, $H_\infty$ is a SAGBI basis
for $R[H_\infty] = R[F]$.

Now suppose that $R[F]$ has a finite SAGBI basis $S$.
Because $H_\infty$ is also a SAGBI basis, for each $s\in S$,
we have an expression
$$\lt(s)=\sum_{j=1}^{M_s} r_{j,s}\lt(H_\infty^{\e_{j,s}}),{~~}r_{j,s}\in R.$$
The finite set $\widetilde{H}$ of those elements of $H_\infty$ for which the
corresponding coordinate of some $\Ve_{j,s}$ is non-zero is 
a SAGBI basis as well since $R[\Lt \widetilde{H}]=R[\Lt S]=R[\Lt R[F]]$.
The set $\widetilde{H}$ must be a subset of the set $H=H_{N_0}$ produced
at the end of some finite number $N_0$ of loops, so that $H_{N_0}$
is also a SAGBI basis for $R[F]$, and by Theorem \ref{314}, we know that
the algorithm will terminate after the next loop.

It remains to show that $H_{N_0}$ is finite.
Any loop that begins with a finite input set (as does the very first
loop, by assumption on $F$) will create a finite associated variable
set $Y$.  Then the Hilbert Basis Theorem applies to $R[Y]$ to prove that
we can choose the generating set $\cal P$ to be finite as well.
Hence, the output of that pass of the loop must be finite.
Thus, beginning with a finite set $F$, Algorithm \ref{316}
completes a strictly finite number of loops, each of which yields finite
output, and we conclude that $H_{N_0}$ is indeed a finite SAGBI
basis for $R[F]$.
\end{pf}

\begin{ex}\label{3999}
In this example we will compute a SAGBI
basis for $\Z[F]\subseteq \Z[x,y]$ where 
$$F=\{
f_1=4x^2y^2+2xy^3+3xy,{~~}f_2=2x^2+xy,{~~}f_3=2y^2\}.$$
We use the term order degree lex with $x>y$.

Set $H=F$.
It is evident that the ideal of relations $I(\Lt H)=I(4x^2y^2,2x^2,2y^2)$
in $\Z[Y_1,Y_2,Y_3]$ is generated by $P(Y)=Y_1-Y_2Y_3$.
The polynomial $P(H) = 3xy$ cannot be s-reduced via $H$,
so that $red{\cal P}=\{3xy\}$.
This forces a second pass through the WHILE loop with
an additional member $f_4=3xy\in H$.\footnote
{The reader may notice that $f_4$ 
is actually an s-reductum of $f_1$
via $\{f_2,f_3\}$ and that we could therefore
have replaced $f_1$ by $f_4$ before beginning the
computations at all.
Such inter-reduction and replacement may well 
make the algorithm more efficient.
However, a serious analysis of 
improvements is outside the scope of
this exposition.  The present example is intended merely to illustrate
our basic algorithm, so we will avoid
introducing any extra techniques at this juncture, tempting and
helpful though it may be.}

On the second pass through the WHILE loop, we calculate generators for
the new ideal $I(\Lt H)\subseteq \Z[Y_1,Y_2,Y_3,Y_4]$, obtaining the set
$${\cal P} =\{P_1=Y_1-Y_2Y_3,{~~}P_2=9Y_1-4Y_4^2,{~~}P_3=
9Y_2Y_3 - 4Y_4^2\}.$$
One can check that $P_j(H)\redH 0$ for $j=1,2,3$.
Thus, the set $red\cal P$ of non-zero s-reducta of
$\cal P$ is empty, terminating the algorithm.
Our SAGBI basis is $\{4x^2y^2+2xy^3+3xy,2x^2+xy,2y^2,3xy\}$.
\end{ex}

\section{SAGBI-\gr Bases in $\r$}
We next address the topic of SAGBI-\gr bases
in $\r$
and begin by defining the primary objects of study.
Then we present an algorithm for their construction.
As always, $A$ is an $R$-\sa of $\r$.
\begin{df}\label{401}
Let $I\subseteq A$ be an ideal of $A$.  A subset $G\subseteq I$
is a SAGBI-\gr basis (SG-basis) for $I$ if 
$\Lt G$ generates $\<\Lt I\>$ in $R[\Lt A]$.
\end{df}
Recall that in ordinary \gr basis theory
every ideal is assured to have a finite \gr basis, due to the
Hilbert Basis Theorem.
By the same reasoning, we may draw this
conclusion about SG-bases for ideals of $A$ 
provided that $A$ has a finite
SAGBI basis.

We continue by describing an appropriate
reduction theory for the current context.

\begin{df}\label{404}
Let $G\subseteq A, h\in A$.  We say that $h$ si-reduces
to $h'$ via $G$ in one step, written $h \sig h'$, if there exist a
non-zero term $cX^\val$ of $h$ and elements
$g_1, \ldots g_M\in G$ and $a_1, \ldots, a_M\in A$
for which the following hold:
\begin{enumerate}
\item $X^\val = \lp(a_ig_i)$ for each $i$.
\item $cX^\val = \sum_{i=1}^M \lt(a_ig_i)$.
\item $h'= h - \sum_{i=1}^M a_ig_i$.
\end{enumerate}
We say that $h$ si-reduces to $h'$ via $G$ and again
write $h\sig h'$
if there is a chain of one-step si-reductions
as above leading from $h$ to $h'$.
If $h'$ cannot be si-reduced via $G$, we call it
a final si-reductum of $h$.
\end{df}
We point out that $h\sig h'$ implies that $h-h'\in \<G\>_A$.
Again, well-ordering of $\T_X$ implies
that every $h\in A$ must have a final si-reductum
via a subset $G$; that is, si-reduction always terminates.

To perform si-reduction, given a term $cX^\val$
of $h$,
we must determine
\begin{enumerate}
\item for each $g\in G$,
whether $X^\val =\lp(g)\lp(a)$ for some $a\in A$, that is,
whether $X^\val\in \<\lp(g)\>_{\sLp A}$, and
\item whether $c$ can be expressed as an $\Lc A$-linear combination
of the appropriate $\lc(g)$'s. (This is equivalent to Condition 2 
of Definition \ref{404} under the homogeneity of Condition 1.)
\end{enumerate}
Given a SAGBI basis $F$ for $A$,
answering the monoid-ideal membership question posed first
amounts to searching for solutions $\vec{\eta}
\in \oplus_F\NN$ to the equation
$$X^\val = \lp(g)\lp(F^\veta)$$
for each $g\in G$, which may be converted to an inhomogeneous
linear diophantine system in its exponents.
We can then check the desired property for the coefficient $c$,
by our assumption that ideal membership in $R$ can be determined.

The proofs of the next result and its corollaries again proceed in the
standard way.
\begin{prop}\label{405}
The following are equivalent for a subset $G$ of an ideal $I\subseteq A$:
\begin{enumerate}
\item $G$ is an SG-basis for $I$.
\item For every $h\in I$, every final si-reductum of $h$ via $G$ is 0.
\item Every $h\in I$ has what is called an 
SG-representation with respect to $G$, that is, a representation
$$h=\sum_{i=1}^M
a_i g_i,{~} a_i\in A,g_i\in G$$
such that $\max_i \lp(a_ig_i)=\lp(h)$.
\end{enumerate}
\end{prop}
\begin{cor}
An SG-basis for $I$ generates $I$ as an ideal of $A$.
\end{cor}
\begin{cor}
Suppose that $G$ is an SG-basis for $I\subseteq A$.
Then $a\in A$ belongs to $I$ $\iff$ $a\sig 0$.
\end{cor}

We introduce some basic terminology.

\begin{df}\label{4211}
For a vector $\Va\in \oplus_G A$ whose coordinates are denoted by $a_g$,
we write $\vlt(\Va)$ for the vector in $\oplus_G\Lt A$
whose $g$-th coordinate is $\lt(a_g)$.
\end{df}

\begin{df}\label{421}
$\lts^*(G)=\{\Va\in \oplus_G A:\vlt(\Va)\in
\syz^*(\Lt G)\subseteq \oplus_G R[\Lt A]\}$.
An element of $\lts^*(G)$ is called a homogeneous $A$-lt-syzygy for $G$.
\end{df}

\begin{df}\label{4215}
We call ${\cal Q}\subseteq\lts^*(G)$ an lt-generating set
for $\lts^*(G)$ if $\{\vlt(\VQ):\VQ\in{\cal Q}\}$ is a
generating set for $\syz^*(\Lt G)$.
\end{df}

For the remainder of this section we assume that
$A$ has a finite SAGBI basis, and that
$G
=\{g_1,\ldots,g_M\}\subseteq A$ is finite
as well; this assures computability.
Given an lt-generating set ${\cal Q}$ and writing its elements as 
$\VQ_j=(q_{j,1},\ldots,q_{j,M})$, we shall see
that the polynomials $\sum_{i=1}^M q_{j,i}g_i$
take the place of S-polynomials in our present setting.

\begin{th}\label{423}
Let $G=\{g_1,\ldots,g_M\}\subseteq A$; let
${\cal Q}$ be an 
lt-generating set for $\lts^*(G)$.
Then $G$ is an SG-basis for $\<G\>_A$ $\iff$
for each $\VQ_j=(q_{j,1},\ldots,q_{j,M})\in {\cal Q}$,
we have $\sum_{i=1}^M q_{j,i}g_i\sig 0$.
\end{th}
\begin{pf}
$\Longrightarrow$: The result is a direct consequence of
Proposition \ref{405}.

$\Longleftarrow$:
Let $h\in \<G\>_A$; write $h=\sum_{i=1}^m a_ig_i$ such that
the height $t_0=\max_i\lp(a_ig_i)$ of this representation is
minimal with respect to all such representations for $h$.
Now $\lp(h)\leq t_0$;
suppose that $\lp(h)< t_0$.
Without loss of generality, assume that our representation
is written such that $a_1g_1,\ldots , a_{M_0}g_{M_0}$ contribute
to the height, in the sense of Definition \ref{315}.
Setting $\Vapr=(a_1,\ldots,a_{M_0},0,\ldots,0)$, we see that
$\vlt(\Vapr)\in\syz^*(\Lt G)$.
Thus, there exist
$b_1,\ldots,b_N\in A$ and $\VQ_1,\ldots,\VQ_N\in{\cal Q}$
such that
$\vlt(\Vapr)=\sum_{j=1}^N \lt(b_j) \vlt(\VQ_j)$;
also, we may assume that
$\deg(\lt(b_j)\vlt(\VQ_j))
=\deg(\vlt(\Vapr))=t_0$ for all $j$
by homogeneity of the syzygies involved.
Furthermore, the elements $b_j$ and $\VQ_j$ may be chosen
so that
the expression $\sum_{j=1}^N \lt(b_j)\lt(q_{j,i})$ is
homogeneous in $\r$
for all $i$ since every non-zero $\lt(b_j)\lt(q_{j,i})\lt(g_i)=
\deg(\lt(b_j)\vlt(\VQ_j))=t_0$.

Now
\begin{eqnarray}
\nonumber h&=&\sum_{i=1}^M a_ig_i -\sum_{i=1}^M(\sum_{j=1}^N b_jq_{j,i})g_i
+\sum_{j=1}^N b_j(\sum_{i=1}^Mq_{j,i}g_i)\\
&=& \sum_{i=1}^M(a_i-\sum_{j=1}^N b_j q_{j,i})g_i
+\sum_{j=1}^N b_j (\sum_{i=1}^M p_{j,i}g_i)\label{5009}
\end{eqnarray}
where $\sum_{i=1}^M p_{j,i}g_i$ is an SG-representation for
$\sum_{i=1}^M q_{j,i}g_i$, which exists since we have supposed that
every $\sum_{i=1}^M q_{j,i}g_i \sig 0$.
Furthermore, if we define $t_j=\het(\sum_{i=1}^M p_{j,i}g_i)$,
then we have 
$$t_j=
\lp(\sum_{i=1}^M q_{j,i}g_i)<\max_i \lp(q_{j,i}g_i) \forall j,$$
where the inequality
holds because $\VQ_j\in\lts^*(G)$.

We proceed to show that the representation for $h$
in Equation (\ref{5009})
has lesser height than our original representation.
The height of the first sum (indexed by $i$) 
is $\max_i \lp[(a_i-\sum_{j=1}^N b_jq_{j,i})g_i]$.
For $i\leq M_0$, we know that $\lt(a_i)=\lt(\sum_{j=1}^N b_j q_{j,i})$
by homogeneity of $\sum_{j=1}^N b_j q_{j,i}$ in $\r$;
therefore, due to cancellation of the highest terms, 
$\lp[(a_i-\sum_{j=1}^N b_j q_{j,i})g_i]<\lp(a_ig_i)=t_0$, our original height.
For $i>M_0$, we recognize that $\lp[(a_i-\sum_{j=1}^N b_jq_{j,i})g_i] \leq \max\{\lp(a_ig_i),
\lp(\sum_{j=1}^N b_j q_{j,i}g_i)\}$, for we assume that
the expression $a_i-\sum_{j=1}^N b_j q_{j,i}$ represents a
simplified polynomial in $\r$.
Yet $i>M_0$ implies that $\lp(a_ig_i)<t_0$ and that
$\sum_{j=1}^N \lt(b_j)\lt(q_{j,i})=0$, which in turn implies that
$\lp(\sum_{j=1}^N b_jq_{j,i}g_i)<\max_j \lp(b_jq_{j,i}g_i)=
\deg(\lt(b_j)\vlt(\VQ_j))=t_0$.
Thus, the height of the first sum in Equation (\ref{5009}) must be less than
the original height since for all $i$, $\lp[(a_i-\sum_{j=1}^N b_j q_{j,i})g_i]
<t_0$.

Now for the second sum, we have the following:
\begin{eqnarray*}
\het(\sum_{j=1}^N b_j\sum_{i=1}^M p_{j,i}g_i)
&\leq& \max_{i,j}\lp(b_jp_{j,i}g_i)=\max_j[\lp(b_j)\cdot t_j]\\
&<&\max_{i,j} \lp(b_jq_{j,i}g_i)=\deg(\vlt(\Vapr))=t_0
\end{eqnarray*}

Hence, Equation (\ref{5009}) does provide
 a new representation for $h\in \<G\>_A$ 
having smaller height than our assumed minimum.
Therefore, $\lp(h)=t_0$, the minimum possible height, proving that
 $G$ is an
SG-basis for $\<G\>_A$.
\end{pf}

We next describe how an lt-generating set for $\lts^*(G)$
may be computed (when $G=\{g_1,\ldots,g_M\}$ is finite).
The method is based on
the following result, whose proof is straight-forward.

\begin{prop}\label{480}
Let $\pi:{\cal R}\longrightarrow {\cal S}$ be a ring epimorphism.
Let $S'=\{s_1,\ldots,s_M\}\subseteq \cal S$ be given, and choose 
a set $\bar{S'}=\{\bar{s}_1,
\ldots,\bar{s}_M\}$ of pre-images in $\cal R$.
Suppose that $\VP_1,\ldots,\VP_L\in {\cal R}^M$ with
$\VP_j=(p_{j,1},\ldots,p_{j,M})$
are such that 
$$\VP_1,\ldots,\VP_K \mbox{ generate }\syz(\bar{s}_1,
\ldots,\bar{s}_M)\subseteq {\cal R}^M$$
while for the remaining
$\{\VP_{K+1},\ldots,\VP_L\}$,
$$\sum_{i=1}^M p_{K+1,i}\bar{s}_i, \ldots, 
\sum_{i=1}^M p_{L,i}\bar{s}_i
\mbox{ generate }\ker(\pi)\cap\<\bar{s}_1,
\ldots,\bar{s}_M\>\subseteq \cal R.$$
Then
$Syz(s_1,\ldots,s_M)$ is generated by the set
$\{\vpi(\VP_1),\ldots,\vpi(\VP_L)\},$
where we define $\vpi:{\cal R}^M\longrightarrow {\cal S}^M$ via
$\vpi(r_1,\ldots,r_M)=(\pi(r_1),\ldots,\pi(r_M))$
for $r_1,\ldots,r_M\in \cal R$.
\end{prop}

To apply this result in the desired setting, we take ${\cal S}=R[\Lt A]=
R[\Lt F]$ where $F$ is a finite SAGBI basis for $A$, set
${\cal R}=R[Y]$
where $Y$ is a set of variables of the same cardinality as $F$,
and take $\pi$ to be the obvious evaluation map.
Proposition \ref{480} and ordinary \gr basis techniques
then allow us to compute generators for $\syz(\Lt G)$, from which
we may obtain a homogeneous generating set ${\cal P}=
\{\VP_1,\ldots,\VP_N\}$ for $\syz^*(\Lt G)$.
Furthermore, we may assume that for each generator $\VP_j(\Lt F)=
(P_{j,1}(\Lt F),\ldots, P_{j,M}(\Lt F))$, the polynomials $P_{j,i}(\Lt F)$
are homogeneous in $\r$.  Defining 
\[
\widetilde{P_{j,i}}(F)=
\left\{\begin{array}{ll}
P_{j,i}(F)&\mbox{ if $P_{j,i}(\Lt F)\neq 0$}\\
0&\mbox{ otherwise,}
\end{array}
\right.
\]
we see that the set
${\cal Q}=\{(\widetilde{P_{j,1}}(F),\ldots,\widetilde{P_{j,M}}(F)):
j=1,\ldots, N\}$ is an lt-generating set for $\lts^*(G)$,
for $\lt(\widetilde{P_{j,i}}(F))=P_{j,i}(\Lt F)$ for all $i$ and $j$.

Next we present an algorithm for computing SG-bases.  See Algorithm \ref{60000}.

\begin{figure}[hbt]
\begin{center}
\frame{
\begin{picture}(380,245)
\normalsize{
\put(30,220){{\bf INPUT:} A finite set $G\subseteq A$, $F$ a finite SAGBI basis for $A$}
\put(30,195){{\bf OUTPUT:} An SG-basis $H$ for $\<G\>_A$}
\put(30,170){{\bf INITIALIZATION:} $H:=G,{~}oldH:=\emptyset$}
\put(30,145) {{\bf WHILE} $H \neq oldH$ {\bf DO}}
\put(50,120) {Compute an lt-generating set $\cal Q$ for $\lts^*(H)$.}
\put(50,95) {${\cal P}:=\{\sum_{h\in H} q_h h:(q_h)_{h\in H}\in
{\cal Q}\}$}
\put(50,70) {$red{\cal P}:=\{\mbox{final si-reducta via $H$ of each element of
 $\cal P$}\}-\{0\}$}
\put(50,45) {$oldH:=H$}
\put(50,20) {$H:= H \cup red{\cal P}$}
}
\end{picture}
}
\begin{mycaption}{SG-Basis Construction Algorithm}\label{60000}
\end{mycaption}
\end{center}
\end{figure}

\begin{prop}\label{60001}
Algorithm \ref{60000} yields a finite SG-basis for $\<G\>_A$
(when $G$ is finite and $A$ has a finite SAGBI basis).
\end{prop}
\begin{pf}
We first show that the algorithm produces an SG-basis,
then that the resulting basis is finite.

Set $H_\infty = \cup H$ 
over all passes of the WHILE loop.
For each ${\cal Q}$ of each loop, construct a set
${\cal Q}'\subseteq \oplus_{H_\infty}A$ by adding
sufficiently many 0 coordinates to each vector in ${\cal Q}$.
We claim that the set ${\cal Q}_\infty'=\cup{\cal Q}'$
over all passes of the loop is an lt-generating set
for $\lts^*(H_\infty)$.
For choose $\vec{\tau}=(t_i)_{h_i\in H_\infty}\in\syz^*(\Lt H_\infty)$.
Only finitely many coordinates $t_i$ are non-zero, corresponding
to finitely many elements $h_i\in H_\infty$.
These elements all belong to the set $H=H_{N_0}$
produced at the end of some finite number of passes of the
WHILE loop.  Defining $\vec{\tau}_0$ to be the
vector consisting precisely of the non-zero coordinates of
$\vec{\tau}$, we note that $\vec{\tau}_0\in \syz^*(H_{N_0})$ and therefore
belongs to the $R[\Lt A]$-module generated by $\{\vlt(\VQ):
\VQ\in {\cal Q}_{N_0}\}$, where ${\cal Q}_{N_0}$ is the chosen
lt-generating set for $\lts^*(H_{N_0})$.
Consequently, $\vec{\tau}$ belongs to the $R[\Lt A]$-module
generated by $\{\vlt(\VQpr):\VQpr
\in {\cal Q}_\infty'\}$,
proving the claim.

We next show that $H_\infty$ is an SG-basis for $\<G\>_A$.
Choose $\VQpr\in{\cal Q}_\infty'$.  Again, $\VQpr=(q_h')_{h\in H_\infty}$
has only finitely many non-zero coordinates, corresponding to a finite
subset of some $H=H_{N_0}\subseteq H_\infty$.
Clearly, $\sum_{q_h'\neq 0} q_h' h$ si-reduces
to 0 via some subset of $H_\infty$, hence via $H_\infty$
either in the loop in which $\VQpr$ is created or in the next.
Thus, $H_\infty$ satisfies Theorem \ref{423}, proving that it
is indeed an SG-basis for $\<G\>_A$.

Since $A$ has a finite SAGBI basis, we know that there exists a
finite SG-basis $S$ for $\<G\>_A$.
We have shown above that $H_\infty$ is an
SG-basis for $\<G\>_A$; therefore,
it must be that for each $s\in S$ there exist
$h_{1,s},\ldots,h_{M_s,s}\in H_\infty$ and
$a_{1,s},\ldots,a_{M_s,s}\in A$ such that
$$\lt(s)=\sum_{i=1}^{M_s} \lt(a_{i,s}h_{i,s}).$$
The set $\widetilde{H}=\cup_{s\in S} \{h_{1,s},\ldots,h_{M_s,s}\}$
is clearly finite, and it is an SG-basis for $\<G\>_A$
since $\<\Lt \widetilde{H}\>=\<\Lt S\>=\<\Lt G\>\subseteq{R[\Lt A]}$.
Because $\widetilde{H}$ must be a subset of the set $H=H_{N_0}$
produced after some finite number of passes of the WHILE loop, $H_{N_0}$
is also an SG-basis, and the algorithm will terminate
at the next loop.

Finally, we show that $H_{N_0}$ is finite.
Our technique for computing an lt-generating set for $\lts^*(H)$
involves calculating a generating set
for $\syz^*(\Lt H)$; these two sets have the same cardinality,
according to Definition \ref{4215}.
Since $R[\Lt A]$ is noetherian,
we may choose a finite generating set for
$\syz^*(\Lt H)$
when the input set $H$
for the loop is finite.
Therefore, ${\cal P}$ and consequently the output of such a loop
are finite.
Then since $H_{N_0}$ is the result of a finite number of passes
of the loop, beginning with finite input $G$, it is a finite
SG-basis for $\<G\>_A$.
\end{pf}

The example below demonstrates how to compute an SG-basis.
\begin{ex}\label{4809}
As in Example \ref{3999},
we take $A={\bf Z}[F]\subseteq {\bf Z}[x,y]$ where 
$$F=\{ f_1=2x^2+xy,{~~}f_2=2y^2,{~~}f_3=3xy\},$$
and let $G\subseteq{\bf Z}[F]$ be given by
$$G=\{g_1=4x^2y^2+2xy^3,{~~}g_2=18x^2y^4\}.$$
We will again use the term order degree lex with $x>y$, with respect
to which we have found that $F$ is a SAGBI basis for $A$.

We begin by setting $H=G$.
Applying the technique described after Proposition \ref{480},
we calculate\footnote{
Some of the intermediate computations
were performed using the Mathematica sub-package
GroebnerZ.  See \cite{NG}.}
 an lt-generating set ${\cal Q}=
\{(f_3^2,-f_1),(9f_2,-4)\}$ for $\lts^*(H)$;
we obtain the associated set ${\cal P}=\{0,36xy^5\}$.
We easily see that $red{\cal P}=\{36xy^5\}$
since this element cannot be si-reduced via $H$.
Therefore, we define $$g_3=36xy^5$$
and conduct a second
pass of the WHILE loop.
This time, we construct
$${\cal Q}=\{
(f_3^2,-f_1,0),(3f_2f_3,0,-f_1),(0,3f_2,-f_3),(-9f_2,4,0)\}.$$
This again yields ${\cal P}=\{0,36xy^5\}$, so
clearly, $red{\cal P}=\emptyset$ now,
and the stopping criterion $H=oldH$ is satisfied.
We have that $$\{4x^2y^2+2xy^3,18x^2y^4,36xy^5\}$$
is an SG-basis for $\<G\>_A$.
\end{ex}

\section{$A$-syzygies}

To conclude, we will present a method for calculating a set of generators
for $\syz_A(H)$ given a finite subset $H$ of an $R$-\sa $A\subseteq\r$,
where we again assume that $A$ has a finite SAGBI basis.
Our technique is based on the following theorem:

\begin{th}\label{900}
Let $G=\{g_1,\ldots,g_M\}\subseteq A$ be a finite
SG-basis for $\<G\>_A$.  Let ${\cal Q}=\{\VQ_1,\ldots,\VQ_N\}$
be an lt-generating set for $\lts^*(G)$,
and write each $\VQ_j=(q_{j,1},\ldots,q_{j,M})$.
For each $j$, let $\sum_{i=1}^M p_{j,i}g_i$
be an SG-representation for $\sum_{i=1}^M q_{j,i}g_i$.
Then $\syz_A(G)$ is generated as an $A$-module
by the vectors 
$$\VP_j=(q_{j,1}-p_{j,1},\ldots,q_{j,M}-p_{j,M}),{~~~}
j=1,\ldots, N.$$
\end{th}
\begin{pf}
Let
$\cal M$ represent the $A$-submodule of $\syz_A(G)$
generated by the set $\{\VP_1,\ldots,\VP_N\}$, and
suppose that the conclusion of the theorem is false.
Then we can choose $\Vh=(h_1,\ldots, h_M)\in\syz_A(G)-{\cal M}$,
such that $t_0=\het(\sum_{i=1}^M h_ig_i)$ as defined in 
Definition \ref{315}
is minimal among such elements of $\syz_A(G)$.
Without loss of generality, we assume that precisely $h_1,\ldots, h_{M_0}$
contribute to the height of this expression.
This implies that $\sum_{i=1}^{M_0} \lt(h_i)\lt(g_i)=0$, i.e.,
that $\vlt(\Vhpr)\in\syz^*(\Lt G)\subseteq R[\Lt A]$
where $\Vhpr=(h_1,\ldots,h_{M_0},0,\ldots,0)$.
Therefore, we can write
$$\vlt(\Vhpr)=\sum_{j=1}^N \lt(b_j)\vlt(\VQ_j)$$
where $\deg[\lt(b_j)\vlt(\VQ_j)]=\deg(\vlt(\Vhpr))=t_0$
for all $j$ such that $b_j\neq 0$.
Also, as we saw in the proof of Theorem \ref{423},
we may assume that the expression $\sum_{j=1}^N \lt(b_j)\lt(q_{j,i})$
is homogeneous in $\r$ for all $i$.

Now we consider the element $\Vs=\Vh -\sum_{j=1}^N b_j \VP_j\in
\syz_A(G)-{\cal M}$.  We claim that $\deg(\Vs)<t_0$,
By definition, $\deg(\Vs)=\max_i \lp(s_ig_i)$ where $s_i$ is
the $i$-th coordinate of $\Vs$; in particular, $s_i$ is the
simplified form of $h_i-\sum_{j=1}^N b_j p_{j,i}$.
For $i\leq M_0$, $\lt(\sum_{j=1}^N b_j p_{j,i})=\lt(\sum_{j=1}^N b_j q_{j,i})
=\lt(h_i)$; whence, cancellation of the highest
terms yields $\lp(s_ig_i)<\lp(h_ig_i)=t_0$.
For $i>M_0$, $\lp(s_ig_i)\leq \max\{\lp(h_ig_i),
\lp(\sum_{j=1}^N b_jp_{j,i}g_i)\}$.
By assumption, $\lp(h_ig_i)<t_0$, and
$$\lp(\sum_{j=1}^N b_jp_{j,i}g_i)=\lp(\sum_{j=1}^N b_jq_{j,i}g_i)
<\max_j \lp(b_jq_{j,i}g_i)=t_0$$
where the inequality holds because $\sum_{j=1}^N \lt(b_j q_{j,i})=0$
for $i>M_0$ and the final equality holds due to our assumption
that $\deg(\lt(b_j)\vlt(\VQ_j))=\deg(\vlt(\Vhpr))=t_0$.
Hence, $\lp(s_ig_i)<t_0$ for all $i$, and we indeed have
$\deg(\Vs)<t_0$,
which contradicts our assumption of the minimality of $t_0$
for elements of $\syz_A(G)-{\cal M}$.
Therefore, this difference is empty, and $\syz_A(G)={\cal M}$.
\end{pf}

We are now prepared to
compute the generators for the $A$-syzygy module $\syz_A(H)$
of an arbitrary finite subset $H\subseteq A$.
We briefly outline the standard technique, which is described in 
greater detail in
such references as \cite{AL}.
Specifically,
we
compute an SG-basis $G$ for $\<H\>_A$
and then produce matrices $\cal W$ and $\cal U$
with entries in $A$
such that $H={\cal W}G$ and $G={\cal U}H$,
where we now view $G$ and $H$ as column vectors.
The module $\syz_A(H)$ is then generated
by the vectors $\VP_j\cal U$
together with the row vectors of ${\cal I}-{\cal WU}$,
where ${\cal I}$ is the identity matrix of the appropriate size.

\begin{ex}
Again, we take $A={\bf Z}[F]\subseteq{\bf Z}[x,y]$
where 
$$F=
\{f_1=2x^2+xy,{~~}f_2=2y^2,{~~}f_3=3xy\}$$
 is a SAGBI
basis for $A$ with respect to our 
term order, degree lex with $x>y$.
Let 
$$H=\{h_1=4x^2y^2+2xy^3,{~~}h_2=10x^2y^4+4xy^5,{~~}h_3=36xy^5\}\subseteq A.$$
It is apparent that the set 
$$G=\{g_1=4x^2y^2+2xy^3,{~~}g_2=18x^2y^4,{~~}g_3=36xy^5\}$$
of Example \ref{4809} is an SG-basis for $\<H\>_A$, for we observe 
that $h_1=g_1,$ $h_2=g_2-f_2g_1,$ and $h_3=g_3$.
Thus, we have the change-of-basis matrices
\[
{\cal W}=\left[
\begin{array}{ccc}
1&0&0\\
-f_2&1&0\\
0&0&1
\end{array}
\right]\mbox{ and }{~}
{\cal U}=\left[
\begin{array}{ccc}
1&0&0\\
f_2&1&0\\
0&0&1
\end{array}
\right]
\]
described above.
Because ${\cal I}-{\cal WU}$ is the zero-matrix, the only non-trivial
generators for $\syz_A(H)$ are the vectors $\VP_j{\cal U}$,
which we will now compute.

We recall the lt-generating set 
\begin{eqnarray*}
{\cal Q}&=\{& \VQ_1 =(f_3^2,-f_1,0),\\
&&\VQ_2=(3f_2f_3,0,-f_1),\\
&&\VQ_3=(0,3f_2,-f_3),\\
&&\VQ_4=(-9f_2,4,0){~}\}
\end{eqnarray*}
for $\lts^*(G)$ as described at the end of Example \ref{4809}.
For the first three of these vectors, 
the polynomials $\sum_{i=1}^3 q_{j,i}g_i =0$;
thus, $\VP_j=\VQ_j$ for $j=1,2,3$.
However, $\VQ_4$ gives us the expression $-9f_2g_1+4g_2=-36xy^5=-g_3$,
which yields $\VP_4 =(-9f_2,4,1)$.
We conclude that
\begin{eqnarray*}
\VP_1{\cal U}&=&(f_3^2-f_1f_2,-f_1,0)\\
\VP_2{\cal U}&=&(3f_2f_3,0,-f_1)\\
\VP_3{\cal U}&=&(3f_2^2,3f_2,-f_3)\\
\mbox{and}{~~}\VP_4{\cal U}&=&(-5f_2,4,1)
\end{eqnarray*}
generate $\syz_A(H)$ as an $A$-module.
\end{ex}

\section{Acknowledgements}
This paper is derived from my doctoral dissertation,
directed by Dr. William Adams, to whom I am deeply grateful.
Thanks also go to Dr. Philippe Loustaunau for his 
helpful comments.


\begin{thebibliography}{9999}

\bibitem[AL]{AL} W. Adams and P. Loustaunau, {\em An
Introduction to \gr Bases}, American Mathematical Society, Providence, 1994.

\bibitem[Dach]{Dach} T. Dachsel, 
Diplomarbeit, Universit\"at Kaiserslautern,
1990.

\bibitem[KM]{KM} D. Kapur and K. Madlener, {\em A Completion Procedure
for Computing a Canonical Basis for a $k$-Subalgebra}. {\em Computers and
Mathematics}, Springer, New York, 1989, 1-11.

\bibitem[Mil]{Mil} J. L. Miller, {\em Algorithms for Computing
in Subalgebras of Polynomial Algebras over a Ring}, doctoral dissertation,
University of Maryland, 1994.

\bibitem[NG]{NG} G. Nakos and N. Glinos, {\em Computing \gr Bases over
{\bf Z}}, preprint, 1994.

\bibitem[Oll]{Oll} F. Ollivier, {\em Canonical bases: relations with
standard bases, finiteness conditions and application to tame automorphisms}.
\`{a} para\^{i}tre dans les actes de MEGA '90, Castiglioncello,
Birkhauser, 1990.

\bibitem[RS]{RS} L. Robbiano and M. Sweedler, {\em Subalgebra Bases}.
Proc. Commutative Algebra Salvador (W. Burns and A. Simis eds.),
Springer LNM {\bf 1430} (1988), 61-87.

\bibitem[IVR]{IVR} M. Sweedler, {\em Ideal bases and valuation rings},
preprint, 1988.

\end{thebibliography}
\end{document}